\documentclass[a4paper,12pt]{article}

\usepackage[top=3.0cm,bottom=3.0cm,left=2.45cm,right=2.45cm]{geometry}%{geometry} 表示可以自定义页面设置
\usepackage{amsmath,amsthm,mathrsfs,graphicx,amsfonts}%{mathrsfs}用于产生一种数学用的花体字，数学公式宏包
\usepackage{graphicx}
\usepackage{amssymb}
\usepackage{amsmath}
\usepackage{subfigure}
\usepackage{graphicx}
\usepackage{epstopdf}
\usepackage{bibentry}
\usepackage{color}
\usepackage{breqn}
\usepackage{lineno,hyperref}
\usepackage{algorithm}      % 算法浮动体环境
\usepackage{algpseudocode}
\usepackage{indentfirst}
\allowdisplaybreaks[4]

\newcommand{\eqa}{\begin{eqnarray}}
\newcommand{\eeqa}{\end{eqnarray}}
\newcommand{\beq}{\begin{equation}}
\newcommand{\eeq}{\end{equation}}

\newcommand{\bsb}{\begin{subequations}}
\newcommand{\esb}{\end{subequations}}

\newtheorem{theorem}{\bf Theorem}[section]

\newtheorem{problem}[theorem]{Problem}

\newtheorem{definition}[theorem]{Definition}

\begin{document}
\title{Graded discrepancy of graphs and hypergraphs}

\author{Yanling Chen\footnote{Center for Discrete Mathematics, Fuzhou University, Fujian 350003, China. Email: chenyanling\_math@163.com.} ~~ Shuping Huang\footnote{School of Mathematical Sciences, University of Science and Technology of China, Hefei, Anhui 230026, China. Email: hsp@mail.ustc.edu.cn.} ~~ Qinghou Zeng\footnote{Center for Discrete Mathematics, Fuzhou University, Fujian 350003, China. Research supported by National Key R\&D Program of China (Grant No. 2023YFA1010202) and National Natural Science Foundation of China (Grant No. 12371342). Email: zengqh@fzu.edu.cn.}}

\date{}
\maketitle

\begin{abstract}
	This paper studies  the following question of Bollobás and Scott \cite{BS2002}:
	Let $G$ be a graph with $n$ vertices and $p\binom{n}{2}$ edges. What is the smallest $c(p, n)$ such that there is an ordering $v_1, \ldots, v_n$ of the vertices in $G$ with
		$\left|e(\{v_1, \ldots, v_i\})-p\binom{i}{2}\right|\leq c(p, n)$
	for all $i\in \{1,\ldots,n\}$ ?
    We obtain upper and lower bounds for $c(p,n)$ that are both linear in $n$. Furthermore, we generalize the result to $k$-uniform hypergraphs.
\end{abstract}

{\bf Keywords:} graded discrepancy, hypergraph, ordering

\section{Introduction}

Let $G$ be a graph on $n$ vertices and $p=e(G)/\binom{n}{2}$ be the edge density of $G$. Define the \emph{discrepancy} of $G$ to be
\begin{equation*}
	\operatorname{disc}(G)=\max _{S \subseteq V(G)}\left|e(S)-p\binom{|S|}{2}\right|,
\end{equation*}
where $e(S)$ is the number of edges in $G$ induced by $S$.  The discrepancy measures how uniformly the edges of $G$ are distributed among the vertices.  This
important concept appears naturally in various branches of Combinatorics and has been studied by many researchers in recent years.
%The discrepancy is closely related to the theory of quasirandom graphs (see [7])

A fundamental result in this area was established by Erd\H{o}s and Spencer \cite{ES1971}, who proved that every graph $G$ with $n$ vertices has a subset $S\subseteq V(G)$ such that $|e(S)-\frac{1}{2}\binom{|S|}{2}|\geq cn^{3/2}$ for some positive constant $c$, yielding the bound $\operatorname{disc}(G)\ge cn^{3/2}$ for $G$ with $p=1/2$. This result was later extended by Erd\H{o}s, Goldberg, Pach and Spencer \cite{EGP1988} to graphs of arbitrary edge density $p$. A further generalization was obtained by Bollob\'as and Scott \cite{BS2006}, who established the tight bound $\operatorname{disc}(G)\ge c\sqrt{p(1-p)}n^{3/2}$ for any graph $G$ with edge density $p$, where $p(1-p)\geq 1/n$. Remarkably, the random graphs demonstrate that all the above lower bounds are optimal up to constant factors. Bollob\'as and Scott \cite{BS2006} also derived analogous discrepancy bounds for uniform hypergraphs, covering arbitrary edge densities.
%More general results for uniform hypergraphs with arbitrary densities were also proved by Bollob\'as and Scott \cite{BS2006}.
%Beck and S\'os \cite{BS1995}, Bollob\'as and Scott \cite{BS2006}, Chazelle \cite{Cha2000}, Matou\v{s}ek \cite{Mat1999} and S\'os \cite{Sos1983}.
A similar notion is the relative discrepancy of two hypergraphs introduced by Bollob\'as and Scott \cite{BS2011}, measuring the extent to which the edges of the two hypergraphs are uniformly and independently distributed. For more interesting results on discrepancy theory, we refer to  \cite{BS1995,BS2015JCTB,BS2015EJC,Cha2000,MNS2013,Mat1999,RST2023,RT2024,Sos1983}.

The aim of this paper is to study the following so-called \textit{graded discrepancy} problem, introduced by Bollob\'as and Scott (see Problem 9 in \cite{BS2002}).

\begin{problem}[Bollob\'as and Scott \cite{BS2002}]\label{BS-Problem}
	Let $G$ be a graph with $n$ vertices and $p\binom{n}{2}$ edges. What is the smallest $c(p, n)$ such that there is an ordering $v_1, \ldots, v_n$ of the vertices in $G$ with
	\begin{equation*}
		\left|e(\{v_1, \ldots, v_i\})-p\binom{i}{2}\right|\leq c(p, n)
	\end{equation*}
	for all $i\in \{1,\ldots,n\}$ ?
\end{problem}

For an integer $k\ge2$, let $H$ be a $k$-uniform hypergraph with $n$ vertices and $p\binom{n}{k}$ edges. Given an ordering $v_1, \ldots, v_n$ of the vertices in $H$, let $H_i$ denote the subhypergraph induced by $\{v_1,\ldots,v_i\}$ for each $i\in\{1,\ldots,n\}$. This implies that there exists a sequence of induced subhypergraphs with
\begin{equation*}
	H_1 \subsetneq H_2 \subsetneq \ldots \subsetneq H_{n-1}\subsetneq H_n=H.
\end{equation*}
We first define the \textit{local discrepancy} of $H_i$ respect to $H$ by
\begin{equation*}
	\operatorname{grdisc}(H_{i},H)=\left|e(H_i)-p\binom{i}{k}\right|.
\end{equation*}
Then, we further define the \textit{positive local discrepancy }of $H_i$ respect to $H$ by
\begin{equation*}
	\operatorname{grdisc}^{+}(H_{i},H)=\max \left\{e(H_i)-p\binom{i}{k}
	,0\right\},
\end{equation*}
and the \textit{negative local discrepancy} by
\begin{equation*}
	\operatorname{grdisc}^{-}(H_i,H)=\max \left\{p\binom{i}{k}-e(H_i),0\right\}.
\end{equation*}
Hence, $\operatorname{grdisc}(H_i,H)=\max \left\{\operatorname{grdisc}^{-}(H_i,H), \operatorname{grdisc}^{+}(H_i,H)\right\}$ for each $i\in\{1,\ldots,n\}$.  Finally, we give the definition of the graded discrepancy of general hypergraphs. For convenience, we write $[t]:=\{1,\ldots,t\}$ for any integer $t>0$.
\begin{definition}
	For an integer $k\ge2$, let $H$ be a k-uniform hypergraph with $n$ vertices and $p\binom{n}{k}$ edges. The graded discrepancy of $H$ is defined as
	\begin{equation*}
		\operatorname{grdisc}(H)=\min_{\sigma\in S_n}\max_{i\in [n]} \operatorname{grdisc}(H_{\sigma(i)},H),
	\end{equation*}
	where $S_n$ denotes the group of all permutations of $[n]$.
\end{definition}
Therefore, $c(p, n)$ can be written as
\begin{equation*}
	c(p, n)=\max_{e(H)/\binom{n}{k}=p}	\operatorname{grdisc}(H).
\end{equation*}
We solve Problem \ref{BS-Problem} by establishing the following two theorems.
\begin{theorem}\label{thm1}
	Let $G$ be a graph with $n$ vertices and $p\binom{n}{2}$ edges.
	Then the graded discrepancy of $G$ satisfies
	\begin{equation*}
		\operatorname{grdisc}(G) \le
		\begin{cases}
			\displaystyle 	\frac{1+p}{2}  \left(n-1\right), & 0 < p < \frac{1}{3}, \\[10pt]
			\displaystyle 	\max\{p, 1-p\}  \left(n-1\right), & \frac{1}{3} \le p \le \frac{2}{3}, \\[10pt]
			\displaystyle 	(1 - \frac{p}{2}) \left(n-1\right), & \frac{2}{3} < p < 1.
		\end{cases}
	\end{equation*}
	Consequently, for the function $c(p, n)$, we have
	\begin{equation*}
		c(p, n) \le
		\begin{cases}
	\displaystyle 	\frac{1+p}{2}  \left(n-1\right), & 0 < p < \frac{1}{3}, \\[10pt]
	\displaystyle 	\max\{p, 1-p\}  \left(n-1\right), & \frac{1}{3} \le p \le \frac{2}{3}, \\[10pt]
	\displaystyle 	(1 - \frac{p}{2}) \left(n-1\right), & \frac{2}{3} < p < 1.
\end{cases}
	\end{equation*}
\end{theorem}
\begin{theorem}\label{thm2}
	For every  $p \in (0,1)$,  there is a  graph $G$ with $n$ vertices and $p\binom{n}{2}$ edges such that the following holds:
	\begin{equation*}
		\operatorname{grdisc}(G) \ge
		\begin{cases}
\displaystyle  \frac{ 1-p}{2}n+O(1), & 0 < p \le \dfrac{1}{4}, \\[10pt]
	\displaystyle 	\left(\sqrt{p}-p\right)\left(n-1\right)+\sqrt{p}-1, & \dfrac{1}{4} < p \le \dfrac{1}{2}, \\[10pt]
	\displaystyle	\left(\sqrt{1-p}-(1-p)\right)\left(n-1\right)+\sqrt{1-p}-1, & \dfrac{1}{2} \le p < \dfrac{3}{4}, \\[10pt]
\displaystyle  \frac{p}{2}n+O(1), & \dfrac{3}{4} \le p < 1.
\end{cases}
	\end{equation*}
	Consequently, for the function $c(p, n)$, we have
	\begin{equation*}
		c(p,n) \ge
		\begin{cases}
	\displaystyle  \frac{ 1-p}{2}n+O(1), & 0 < p \le \dfrac{1}{4}, \\[10pt]
	\displaystyle 	\left(\sqrt{p}-p\right)\left(n-1\right)+\sqrt{p}-1, & \dfrac{1}{4} < p \le \dfrac{1}{2}, \\[10pt]
	\displaystyle	\left(\sqrt{1-p}-(1-p)\right)\left(n-1\right)+\sqrt{1-p}-1, & \dfrac{1}{2} \le p < \dfrac{3}{4}, \\[10pt]
	\displaystyle  \frac{p}{2}n+O(1), & \dfrac{3}{4} \le p < 1.
\end{cases}
	\end{equation*}
\end{theorem}
We also generalize Theorems \ref{thm1} and \ref{thm2} to the hypergraph setting as follows.
%and obtain a hypergraph version result as follows.
\begin{theorem}\label{thm3}
	For an integer $k\ge2$, let $H$ be a $k$-uniform hypergraph with $n$ vertices and $p\binom{n}{k}$ edges.
	Then the graded discrepancy of $H$ satisfies
	\begin{equation*}
		\operatorname{grdisc}(H) \le
		\begin{cases}
			\displaystyle 	 \frac{1+p}{2}\binom{n-1}{k-1}, & 0 < p < \frac{1}{3}, \\[10pt]
			\displaystyle 	\max\{p, 1-p\}  \binom{n-1}{k-1}, & \frac{1}{3} \le p \le \frac{2}{3}, \\[10pt]
			\displaystyle 	\left(1 - \frac{p}{2}\right) \binom{n-1}{k-1}, & \frac{2}{3} < p < 1.
		\end{cases}
	\end{equation*}
	Consequently, for the function $c(p, n)$, we have
	\begin{equation*}
		c(p, n) \le
		\begin{cases}
	\displaystyle 	 \frac{1+p}{2}\binom{n-1}{k-1}, & 0 < p < \frac{1}{3}, \\[10pt]
	\displaystyle 	\max\{p, 1-p\}  \binom{n-1}{k-1}, & \frac{1}{3} \le p \le \frac{2}{3}, \\[10pt]
	\displaystyle 	\left(1 - \frac{p}{2}\right) \binom{n-1}{k-1}, & \frac{2}{3} < p < 1.
\end{cases}
	\end{equation*}
\end{theorem}
\begin{theorem}\label{thm4}
	For every  $p \in (0,1)$,  there is a $k$-uniform hypergraph $H$ with $n$ vertices and $p\binom{n}{k}$ edges such that the following holds:
	\begin{equation*}
		\operatorname{grdisc}(H) \ge
	\begin{cases}
	\displaystyle \left(\frac{1}{2^{k-1}}-\frac{p}{2}\right)\binom{n-1}{k-1}+O(n^{k-2}), &0<p<\dfrac{1}{2^k}, \\[10pt]
	\displaystyle \left(p^{(k-1)/k}-p-o(p)\right)\binom{n-1}{k-1}, & \dfrac{1}{2^k} < p \le \dfrac{1}{2}, \\[10pt]
	\displaystyle \left((1-p)^{(k-1)/k}-(1-p)-o(1-p)\right)\binom{n-1}{k-1}, & \dfrac{1}{2} \le p < 1-\dfrac{1}{2^k}, \\[10pt]
	\displaystyle\left(\frac{1}{2^{k-1}}-\frac{1-p}{2}\right)\binom{n-1}{k-1}+O(n^{k-2}),  & 1-\dfrac{1}{2^k} \le p < 1.
\end{cases}
	\end{equation*}
	Consequently, for the function $c(p, n)$, we have
	\begin{equation*}
		c(p,n) \ge
	\begin{cases}
	\displaystyle \left(\frac{1}{2^{k-1}}-\frac{p}{2}\right)\binom{n-1}{k-1}+O(n^{k-2}), &0<p<\dfrac{1}{2^k}, \\[10pt]
	\displaystyle \left(p^{(k-1)/k}-p-o(p)\right)\binom{n-1}{k-1}, & \dfrac{1}{2^k} < p \le \dfrac{1}{2}, \\[10pt]
	\displaystyle \left((1-p)^{(k-1)/k}-(1-p)-o(1-p)\right)\binom{n-1}{k-1}, & \dfrac{1}{2} \le p < 1-\dfrac{1}{2^k}, \\[10pt]
	\displaystyle\left(\frac{1}{2^{k-1}}-\frac{1-p}{2}\right)\binom{n-1}{k-1}+O(n^{k-2}),  & 1-\dfrac{1}{2^k} \le p < 1.
\end{cases}
	\end{equation*}
\end{theorem}
\noindent\textbf{Remark}.     In Theorems \ref{thm2} and \ref{thm4}, we assume  that $p = \omega(1/n)$ and $1-p = \omega(1/n)$.

\section{Upper bound of $c(p, n)$}

This section presents the proofs of Theorems \ref{thm1} and  \ref{thm3}, which provide upper bounds for $c(p,n)$. Observing that the graded discrepancies of \(H\) and its complementary graph \(H^c\) must be equal (by taking the same vertex ordering), we immediately obtain the symmetry \(c(p, n) = c(1-p, n)\).
Due to this symmetry, we only need to consider $0<p\leq \frac{1}{2}$. We then develop two algorithms: a general one for any  \( p\in (0,1) \), and a special one for \( p \in (0,\frac{1}{3}) \). Given a hypergraph \( H \), each algorithm generates a sequence of induced subhypergraphs
\begin{equation*}
	H_1 \subsetneq H_2 \subsetneq \ldots \subsetneq H_{n-1} \subsetneq H_n = H .
\end{equation*}
For each sequence, we compute a uniform upper bound of \( \operatorname{grdisc}(H_i, H) \) for all \( i \in [n] \).
Hence, we obtain an upper bound of \( \operatorname{grdisc}(H) \), which  immediately yields an upper bound for $c(p,n)$.

We now present Algorithm 1. Although it is designed for arbitrary $p$, it  achieves its optimal performance in the range $\frac{1}{3}\leq p\leq \frac{1}{2}$.  For a hypergraph \(H\) and a vertex \(v \in V(H)\), let \(H-v\) denote the subhypergraph induced by \(V(H) \setminus \{v\}\), let \(d(v, H)\) be the degree of \(v\) in \(H\), and let \(d(H)\) denote the average degree of \(H\).

~\

\begin{algorithm}[htbp]
	\caption{}
	\label{al1}
	\begin{algorithmic}[1]
		
		\Statex \textbf{Input:}  a $k$-uniform hypergraph  $H$ with $n$ vertices and $p\binom{n}{k}$ edges.
		\Statex \textbf{Output:} a sequence of induced subhypergraphs with
		\begin{equation*}
			H_1 \subsetneq H_2 \subsetneq \ldots \subsetneq H_{n-1}\subsetneq H_n=H.
		\end{equation*}
		
		\State set $i:=n$		
		\While{$i>1$}
		\State  consider  positive local discrepancy $\operatorname{grdisc}^{+}(H_{i},H)$
		
		\If{$\operatorname{grdisc}^{+}(H_{i},H)>0$}
		\State remove a point $v_{i}$ satisfying  $$d(v_i,H_i)=\min _{v\in H_i}\{d(v,H_i) | v \in H_{i+1}, d(v,H_i)\geq d(H_i) \},$$ and then obtain the induced subhypergraph $H_{i-1}=H_i-v$
		\State decrement $i$ by $1$
		\Else
		\State remove a point $v_{i}$ satisfying  $$d(v_i,H_i)=\max _{v\in H_i}\{d(v,H_i) | v \in H_{i+1}, d(v,H_i)\leq d(H_i) \},$$ and then obtain the induced subhypergraph $H_{i-1}=H_i-v$
		\State decrement $i$ by $1$
		\EndIf
		\EndWhile
		
		\State \Return 		\begin{equation*}
			H_1 \subsetneq H_2 \subsetneq \ldots \subsetneq H_{n-1}\subsetneq H_n=H.
		\end{equation*}
		
	\end{algorithmic}
\end{algorithm}

~\

Analyzing the output of Algorithm \ref{al1}, we obtain the following theorem.
\begin{theorem}\label{la1}
	For an integer $k\ge2$, let $H$ be a $k$-uniform hypergraph with $n$ vertices and $p\binom{n}{k}$ edges. Then
	\begin{equation*}
		\operatorname{grdisc}(H)\le \max\{p,1-p\}\binom{n-1}{k-1}.
	\end{equation*}
\end{theorem}
\begin{proof}[{\bf Proof}]
	We only need to prove that the following condition holds at all times during the algorithm's execution:
	\begin{equation}\label{a1}
		\operatorname{grdisc}^+(H_{i},H)\leq p\binom{n-1}{k-1},\quad
		\operatorname{grdisc}^-(H_{i},H)\leq(1-p)\binom{n-1}{k-1},
	\end{equation}
	where $H_i=H_{i+1}-v_{i+1}$  for each $i\in[n-1]$ and $H_n=H$.
	
	First, the condition (\ref{a1}) holds at certain specific steps:

~\

\textbf{Case 1} When $i=n$ and $i\leq k$,   the condition (\ref{a1})  holds trivially.

~\
	
	We now prove that the algorithm preserves the condition (\ref{a1}) at each step. Specifically, we show that if $H_i$ satisfies the condition (\ref{a1}) during the execution, then the constructed hypergraph $H_{i-1}$ also satisfies the condition (\ref{a1}).     Since at most one of $\operatorname{grdisc}^+(H_{i},H)$ and $\operatorname{grdisc}^-(H_{i},H)$ can be nonzero, we consider two possibilities, which are designated as Case 2 and Case 3.
	
	~\
	
	\textbf{Case 2}
	Suppose that $H_i$ satisfies
	$$\operatorname{grdisc}(H_{i},H)=\operatorname{grdisc}^+(H_{i},H)\leq p\binom{n-1}{k-1}, \quad \operatorname{grdisc}^{-}(H_{i},H)=0,\quad i>k.$$
	Then
	\begin{align}\label{i-k-Pos}
		e(H_{i})=p\binom{i}{k}+\operatorname{grdisc}^+(H_{i},H).
	\end{align}
	The average degree of $H_{i}$ is
	\begin{equation*}
		\begin{split}
			d(H_{i})&=\frac{k}{i}\left(p\binom{i}{k}+\operatorname{grdisc}^+(H_{i},H)\right)\\
			&=p\binom{i-1}{k-1}+\frac{k}{i}\operatorname{grdisc}^+(H_{i},H) .
		\end{split}
	\end{equation*}
	We remove a point $v_{i}$ satisfying  $d(H_{i})\leq d(v_{i},H_{i})\leq \binom{i-1}{k-1}$, and obtain the induced subhypergraph $H_{i-1}$.
	It follows that
	\begin{align*}
		e(H_{i-1})&=e(H_{i-1})-d(v_{i},H_{i})\\
		&\le p\binom{i-1}{k}+\frac{i-k}{i}\operatorname{grdisc}^+(H_{i},H).
	\end{align*}
	Then
	\begin{equation*}
		\begin{split}
			\operatorname{grdisc}^+(H_{i-1},H)&=\max \left\{e(H_{i-1})-p\binom{i-1}{k}		,0\right\}\\
			%&\leq\left|p\binom{n-i+1}{k}+\operatorname{grdisc}^+(H_{n-i+1},H)-d(H_{n-i+1})-p\binom{n-i}{k}\right|\\
			&\leq\frac{i-k}{i}\operatorname{grdisc}^+(H_{i},H)\\
			&\leq p\binom{n-1}{k-1}.\\
		\end{split}
	\end{equation*}
	Recall that $e(H_{i-1})=e(H_i)-d(v_{i},H_{i})\ge e(H_{i})-\binom{i-1}{k-1}$. From \eqref{i-k-Pos}, we have
	\begin{equation*}
		\begin{split}
			\operatorname{grdisc}^-(H_{i-1},H)&=\max \left\{p\binom{i-1}{k}-e(H_{i-1}),0\right\}\\
			&\leq \max \left\{p\binom{i-1}{k}-\left(p\binom{i}{k}+\operatorname{grdisc}^+(H_{i},H)-\binom{i-1}{k-1}\right),0\right\}\\
			&=\max\left\{(1-p)\binom{i-1}{k-1}-\operatorname{grdisc}^+(H_{i},H),0\right\}\\
			%&\leq (1-p)\binom{n-i}{k-1}\leq
			&\leq(1-p)\binom{n-1}{k-1}.
		\end{split}
	\end{equation*}

	~\
	
	\textbf{Case 3}  Suppose that $H_i$ satisfies
	$$	\operatorname{grdisc}(H_{i},H)=	\operatorname{grdisc}^-(H_{i},H)\leq(1-p)\binom{n-1}{k-1},\quad \operatorname{grdisc}^{+}(H_{i},H)=0,\quad i>k.$$
	Then
	\begin{align}\label{i-k-Neg}
		e(H_{i})=p\binom{i}{k}-	\operatorname{grdisc}(H_{i},H).
	\end{align}
	The average degree of  $H_{i}$ is
	\begin{align*}
		d(H_{i})&=\frac{k}{i}\left(p\binom{i}{k}-\operatorname{grdisc}^-(H_{i},H)\right)\\
		&=p\binom{i-1}{k-1}-\frac{k}{i}\operatorname{grdisc}^-(H_{i},H).
	\end{align*}
	We remove a point $v_{i}$ satisfying   $0 \leq d(v_{i},H_{i})\leq d(H_{i})\leq p\binom{i-1}{k-1}$, and obtain the induced subhypergraph $H_{i-1}$.
	It follows that
	\begin{align*}
		e(H_{i-1})&=e(H_{i})-d(v_{i},H_{i})\\
		&\ge p\binom{i-1}{k}-\frac{i-k}{i}\operatorname{grdisc}^-(H_{i},H).
	\end{align*}
	This implies that
	\begin{equation*}
		\begin{split}
			\operatorname{grdisc}^-(H_{i-1},H)&=\max \left\{p\binom{i-1}{k}-e(H_{i-1}),0\right\}\\
			%&\leq\left|p\binom{n-i+1}{k}-\operatorname{grdisc}^-(H_{n-i+1},H)-d(H_{n-i+1})-p\binom{n-i}{k}\right|\\
			&\leq \frac{i-k}{i}\operatorname{grdisc}^-(H_{i},H) \\
			&\leq (1-p)\binom{n-1}{k-1}.
		\end{split}
	\end{equation*}
	Note also that $e(H_{i-1})\le e(H_{i})$. This together with \eqref{i-k-Neg} yields that
	\begin{equation*}
		\begin{split}
			\operatorname{grdisc}^+(H_{i-1},H)&=\max \left\{e(H_{i-1})-p\binom{i-1}{k}
			,0\right\}\\
			&\leq\max \left\{p\binom{i}{k}-\operatorname{grdisc}^-(H_{i},H)-p\binom{i-1}{k},0\right\}\\
			&\leq p\binom{i-1}{k-1}\leq p\binom{n-1}{k-1}.
		\end{split}
	\end{equation*}
	This completes the proof of Theorem \ref{la1}.	
\end{proof}

To improve the bound when $p$ is small, we employ Algorithm \ref{al2} for $p < 1/3$, as Algorithm \ref{al1} yields suboptimal estimates in this range.

~\

\begin{algorithm}[htbp]
	\caption{}
	\label{al2}
	\begin{algorithmic}[1]
		
		\Statex \textbf{Input:}  a $k$-uniform hypergraph  $H$ with $n$ vertices and $p\binom{n}{k}$ edges, where $0<p<\frac{1}{3}$.
		\Statex \textbf{Output:} a sequence of induced subhypergraphs with
		\begin{equation*}
			H_1 \subsetneq H_2 \subsetneq \ldots \subsetneq H_{n-1}\subsetneq H_n=H.
		\end{equation*}
		
		\State set $i:=n$		
		\While{$i>1$}
		\State  consider  positive local discrepancy $\operatorname{grdisc}^{+}(H_{i},H)$
		
		\If{$\operatorname{grdisc}^{+}(H_{i},H)\leq  \frac{1-p}{2}\binom{n-1}{k-1}$}
		\State remove a point $v_{i}$ satisfying  $$d(v_i,H_i)=\max _{v\in H_i}\{d(v,H_i) | v \in H_{i+1}, d(v,H_i)\leq d(H_i) \},$$ and then obtain the induced subhypergraph $H_{i-1}=H_i-v$
		\State decrement $i$ by $1$
		\Else
		\State remove a point $v_{i}$ satisfying  $$d(v_i,H_i)=\min _{v\in H_i}\{d(v,H_i) | v \in H_{i+1}, d(v,H_i)\geq d(H_i) \},$$ and then obtain the induced subhypergraph $H_{i-1}=H_i-v$
		\State decrement $i$ by $1$
		\EndIf
		\EndWhile
		
		\State \Return 		\begin{equation*}
			H_1 \subsetneq H_2 \subsetneq \ldots \subsetneq H_{n-1}\subsetneq H_n=H.
		\end{equation*}
		
	\end{algorithmic}
\end{algorithm}

~\

Analyzing the output of Algorithm \ref{al2} leads to the following theorem.
\begin{theorem}\label{la2}
	For an integer $k\ge2$, let $H$ be a $k$-uniform hypergraph with $n$ vertices and $p\binom{n}{k}$ edges, $p\leq\frac{1}{3}$. Then
	\begin{equation*}
		\operatorname{grdisc}(H)\le \frac{1+p}{2}\binom{n-1}{k-1}.
	\end{equation*}
\end{theorem}
\begin{proof}[{\bf Proof}]
	We only need to prove that the following condition holds at all times during the algorithm's execution:
	\begin{equation}\label{a2}
		\operatorname{grdisc}^{+}(H_{i},H)\le \frac{1+p}{2}\binom{n-1}{k-1},\quad 	\operatorname{grdisc}^{-}(H_{i},H)\le \frac{1+p}{2}\binom{n-1}{k-1},
	\end{equation}
	where $H_i=H_{i+1}-v_{i+1}$  for each $i\in[n-1]$ and $H_n=H$.
	
	First, the condition (\ref{a2}) holds at certain specific steps:
	
	~\

\textbf{Case 1} When $i=n$ and $i\leq k$,   the condition (\ref{a2})  holds trivially.

~\

We now prove that the algorithm preserves the condition (\ref{a2}) at each step. Specifically, we show that if $H_i$ satisfies the condition (\ref{a2}) during the execution, then the constructed hypergraph $H_{i-1}$ also satisfies the condition (\ref{a2}). We split the analysis into two cases based on the value of $\operatorname{grdisc}^{+}(H_{i},H)$.

~\

	\textbf{Case 2}
	Suppose that $H_i$ satisfies$$\operatorname{grdisc}^{+}(H_{i},H)\leq  \frac{1-p}{2}\binom{n-1}{k-1},\quad \operatorname{grdisc}^{-}(H_{i},H)\le \frac{1+p}{2}\binom{n-1}{k-1},\quad i>k.$$
	
	~\
	
	Note that at most one of $\operatorname{grdisc}^{+}(H_{i},H)$ and $\operatorname{grdisc}^{-}(H_{i},H)$ is non-zero, so we can break down Case 2 into two subcases.
	
	~\
	
	\textbf{Case 2.1} Suppose that $H_i$ satisfies
	$$\operatorname{grdisc}^{+}(H_{i},H)=0,\quad 0\leq \operatorname{grdisc}^{-}(H_i,H)\le \frac{1+p}{2}\binom{n-1}{k-1},\quad i>k.$$
		Then
	\begin{align*}
		e(H_{i})=p\binom{i}{k}-	\operatorname{grdisc}(H_{i},H).
	\end{align*}
		The average degree of $H_{i}$ is
	\begin{equation*}
		\begin{split}
			d(H_{i})&=\frac{k}{i}\left(p\binom{i}{k}-\operatorname{grdisc}^-(H_{i},H)\right)\\
			&=p\binom{i-1}{k-1}-\frac{k}{i}\operatorname{grdisc}^{-}(H_{i},H).
		\end{split}
	\end{equation*}
Then we remove a point $v_{i}$ satisfying
	$$d(v_i,H_i)\leq d(H_i)= p\binom{i-1}{k-1}-\frac{k}{i}\operatorname{grdisc}^{-}(H_{i},H),$$ and then obtain the induced subhypergraph $H_{i-1}$.
	It follows that
	\begin{equation*}
		\begin{split}
			\operatorname{grdisc}^-(H_{i-1},H)&=\max \left\{p\binom{i-1}{k}-e(H_{i-1}),0\right\}\\
			&\leq p\binom{i-1}{k}-\left(p\binom{i}{k}-\operatorname{grdisc}^{-}(H_i,H)-d(v_i,H_i)\right)\\
			&=\operatorname{grdisc}^{-}(H_i,H)+d(v_i,H_i)-p\binom{i-1}{k-1}\\
			&\leq \frac{i-k}{i} \operatorname{grdisc}^{-}(H_i,H)\leq \frac{1+p}{2}\binom{n-1}{k-1},
		\end{split}
	\end{equation*}
	and
	\begin{equation*}
		\begin{split}
			\operatorname{grdisc}^+(H_{i-1},H)&=\max \left\{e(H_{i-1})-p\binom{i-1}{k}		,0\right\}\\
			%&\leq\left|p\binom{n-i+1}{k}+\operatorname{grdisc}^+(H_{n-i+1},H)-d(H_{n-i+1})-p\binom{n-i}{k}\right|\\
			&=\max\left\{\left(p\binom{i}{k}-\operatorname{grdisc}^{-}(H_i,H)-d(v_i,H_i)\right)-p\binom{i-1}{k}, 0\right\}\\
			&\leq p\binom{i-1}{k-1}\\
			&\leq \frac{1+p}{2}\binom{n-1}{k-1}.
		\end{split}
	\end{equation*}
	
	~\
	
	\textbf{Case 2.2} Suppose that $H_i$ satisfies
	$$0<\operatorname{grdisc}^{+}(H_{i},H)\leq  \frac{1-p}{2}\binom{n-1}{k-1},\quad \operatorname{grdisc}^{-}(H_{i},H)=0,\quad i>k.$$
		Then
	\begin{align*}
		e(H_{i})=p\binom{i}{k}+\operatorname{grdisc}^+(H_{i},H).
	\end{align*}
	The average degree of $H_{i}$ is
	\begin{equation*}
		\begin{split}
			d(H_{i})&=\frac{k}{i}\left(p\binom{i}{k}+\operatorname{grdisc}^+(H_{i},H)\right)\\
			&=p\binom{i-1}{k-1}+\frac{k}{i}\operatorname{grdisc}^+(H_{i},H) .
		\end{split}
	\end{equation*}
	Then we remove a point $v_{i}$ satisfying
	$$d(v_i,H_i)\leq d(H_i)= p\binom{i-1}{k-1}+\frac{k}{i}\operatorname{grdisc}^{+}(H_{i},H),$$ and then obtain the induced subhypergraph $H_{i-1}$.
	It follows that
	\begin{equation*}
		\begin{split}
			\operatorname{grdisc}^-(H_{i-1},H)&=\max \left\{p\binom{i-1}{k}-e(H_{i-1}),0\right\}\\
			&=\max\left\{ p\binom{i-1}{k}-\left(p\binom{i}{k}+\operatorname{grdisc}^{+}(H_i,H)-d(v_i,H_i)\right), 0\right\}\\
			&=\max\{d(v_i,H_i)-p\binom{i-1}{k-1}-\operatorname{grdisc}^{+}(H_i,H), 0\}\\\
			&=0,
		\end{split}
	\end{equation*}
	and
	\begin{equation*}
		\begin{split}
			\operatorname{grdisc}^+(H_{i-1},H)&=\max \left\{e(H_{i-1})-p\binom{i-1}{k}		,0\right\}\\
			%&\leq\left|p\binom{n-i+1}{k}+\operatorname{grdisc}^+(H_{n-i+1},H)-d(H_{n-i+1})-p\binom{n-i}{k}\right|\\
			&=\max\left\{\left(p\binom{i}{k}+\operatorname{grdisc}^{+}(H_i,H)-d(v_i,H_i)\right)-p\binom{i-1}{k}, 0\right\}\\
			&\leq \frac{1+p}{2}\binom{n-1}{k-1}.
		\end{split}
	\end{equation*}
	
	~\

	Finally, we verify that if the hypergraph $H_i$ satisfies
	$$\operatorname{grdisc}^{+}(H_{i},H)> \frac{1-p}{2}\binom{n-1}{k-1},$$
	then the hypergraph $H_{i-1}$ still satisfies the desired condition (\ref{a2}).
	
	~\
	
	\textbf{Case 3} Suppose that $H_i$ satisfies $$\frac{1-p}{2}\binom{n-1}{k-1}<\operatorname{grdisc}^{+}(H_{i},H)\leq \frac{1+p}{2}\binom{n-1}{k-1},\quad \operatorname{grdisc}^{-}(H_{i},H)=0,\quad i>k.$$
		Then
	\begin{align*}
		e(H_{i})=p\binom{i}{k}+\operatorname{grdisc}^+(H_{i},H).
	\end{align*}
	The average degree of $H_{i}$ is
	\begin{equation*}
		\begin{split}
			d(H_{i})&=\frac{k}{i}\left(p\binom{i}{k}+\operatorname{grdisc}^+(H_{i},H)\right)\\
			&=p\binom{i-1}{k-1}+\frac{k}{i}\operatorname{grdisc}^+(H_{i},H) .
		\end{split}
	\end{equation*}
	Then we remove a point $v_{i}$ satisfying
	$$d(v_i,H_i)\geq d(H_i)= p\binom{i-1}{k-1}+\frac{k}{i}\operatorname{grdisc}^{+}(H_{i},H),$$ and then obtain the induced subhypergraph $H_{i-1}$.
	It follows that
	\begin{equation*}
		\begin{split}
			\operatorname{grdisc}^-(H_{i-1},H)&=\max \left\{p\binom{i-1}{k}-e(H_{i-1}),0\right\}\\
			&=\max\left\{ p\binom{i-1}{k}-\left(p\binom{i}{k}+\operatorname{grdisc}^{+}(H_i,H)-d(v_i,H_i)\right), 0\right\}\\
			&=\max\{d(v_i,H_i)-p\binom{i-1}{k-1}-\operatorname{grdisc}^{+}(H_i,H), 0\}\\\
			&\leq\max\{(1-p)\binom{i-1}{k-1}-\operatorname{grdisc}^{+}(H_i,H), 0\}\\\
			&\leq \frac{1-p}{2}\binom{n-1}{k-1} \leq \frac{1+p}{2}\binom{n-1}{k-1},
		\end{split}
	\end{equation*}
	\begin{equation*}
		\begin{split}
			\operatorname{grdisc}^+(H_{i-1},H)&=\max \left\{e(H_{i-1})-p\binom{i-1}{k}		,0\right\}\\
			%&\leq\left|p\binom{n-i+1}{k}+\operatorname{grdisc}^+(H_{n-i+1},H)-d(H_{n-i+1})-p\binom{n-i}{k}\right|\\
			&=\max\left\{\left(p\binom{i}{k}+\operatorname{grdisc}^{+}(H_i,H)-d(v_i,H_i)\right)-p\binom{i-1}{k}, 0\right\}\\
			&\leq \frac{i-k}{i} \operatorname{grdisc}^{+}(H_i,H)\leq \frac{1+p}{2}\binom{n-1}{k-1}.
		\end{split}
	\end{equation*}
	This completes the proof of Theorem \ref{la2}.
\end{proof}

\section{Lower bound of $c(p, n)$}
In this section, we prove Theorems \ref{thm2} and \ref{thm4} by constructing extremal  hypergraphs, thereby obtaining lower bounds for $c(p,n)$. Since $c(p, n)=c(1-p, n)$, we continue (as in the previous section) to work under the assumption that $0<p \le 1/2$.     We begin by constructing a family of extremal graphs for small values of $p$.
\begin{theorem}\label{31}
	Let \( H \) be a \( k \)-uniform hypergraph on \( n \) vertices. Its vertex set is the union of two disjoint sets \( A \cup B \), where
	$$ |A| = x ,\quad  |B| = n-x,\quad 0<<x<\frac{n}{2}. $$
	 The edge set consists of three types of edges:
	\begin{itemize}
		\item  	All \( k \)-element subsets of \( A \);
		\item   All edges formed by taking one vertex from \( A \) and \( k-1 \) vertices from \( B \);
		\item  $O(n^{k-2})$ \( k \)-element subsets of \( B \).
	\end{itemize}
	Then the number of edges in \( H \) is given by
	\begin{equation*}
		e(H)=\binom{x}{k} + x \cdot \binom{n-x}{k-1}+O(n^{k-2}).
	\end{equation*}
	The edge density of $H$ is given by
	\begin{equation*}
		p(H)=\frac{\binom{x}{k} + x \cdot \binom{n-x}{k-1}+O(n^{k-2})}{\binom{n}{k}},
	\end{equation*}
	and
\begin{align*}
	\operatorname{grdisc}(H)\geq  (\frac{1}{2^{k-1}}-\frac{p}{2})\binom{n-1}{k-1}+O(n^{k-2}).
\end{align*}
\end{theorem}
\begin{proof}[{\bf Proof of Theorem \ref{31}}]
	Since \( k \)-uniform hypergraph  $H$ contains only two types of vertices, the optimal algorithm is as follows. This algorithm can be decomposed into $j$ steps. At each step, we select an appropriate value $k_j$($=O(1)$) from a predefined family, remove $k_j$ vertices from set $B$, and then remove one vertex from set $A$. This process is repeated until the algorithm terminates, so we can get
	a sequence of induced subhypergraphs with
	\begin{equation*}
		H_1 \subsetneq H_2 \subsetneq \ldots \subsetneq H_{n-1}\subsetneq H_n=H.
	\end{equation*}
	Denote by $A_i$ the restrictions of the sets $A$ to the subhypergraph $H_i$. We define \( d_A(H_{i}) \) as the degree  of any vertex belonging to set \( A_i \) in the subhypergraph \( H_i \).
	
	Note that
	$$d_A(H_{n-k_1})= \binom{x-1}{k-1} +  \binom{n-x-k_1}{k-1}, $$
	and
	$$ e(H_{n-k_1-1})=p\binom{n-k_1}{k}+\operatorname{grdisc}^{+}(H_{n-k_1},H)-d_A(H_{n-k_1}).$$
	It follows that
		\begin{equation*}
		\begin{split}
			\operatorname{grdisc}^-(H_{n-k_1-1},H)&=\max \left\{p\binom{n-k_1-1}{k}-e(H_{n-k_1-1}),0\right\}\\
		&=d_A(H_{n-k_1})-p\binom{n-k_1-1}{k-1}-\operatorname{grdisc}^{+}(H_{n-k_1},H).
	\end{split}
	\end{equation*}
	Then
	\begin{align*}
			\operatorname{grdisc}(H)&\geq 	\max\left\{\operatorname{grdisc}(H_{n-k_1},H), \operatorname{grdisc}(H_{n-k-1},H)\right\}\\
			&=\max\left\{\operatorname{grdisc}^+(H_{n-k_1},H), \operatorname{grdisc}^{-}(H_{n-k_1-1},H)\right\}  \\		
			&\geq \frac{1}{2}\left(d_A(H_{n-k_1})-p\binom{n-k_1-1}{k-1}\right)\\
			&\geq \frac{1}{2}\left(\binom{x-1}{k-1} +  \binom{n-x-k_1}{k-1}- p\binom{n-k_1-1}{k-1}\right)\\
		&\geq  \binom{\lfloor \frac{1}{2}(n-k_1-1) \rfloor}{k-1}- \frac{p}{2}\binom{n-k_1-1}{k-1},
	\end{align*}
where the last inequality follows from the convexity of $\binom{x}{k-1}$. Since we also have $k_1=O(1)$, it follows that
		$$\operatorname{grdisc}(H)\geq  (\frac{1}{2^{k-1}}-\frac{p}{2})\binom{n-1}{k-1}+O(n^{k-2}).$$
		Thus, the desired result holds.
		\end{proof}
Noting that the construction above exists for all \(0 < p < \frac{1}{2^k}\), we immediately obtain the following theorem.

\begin{theorem}
	There is a	$k$-uniform hypergraph with $n$ vertices and $p\binom{n}{k}$ edges, where $p\in (0,\frac{1}{2^k})$,  that satisfies
	$$		\operatorname{grdisc}(H)\geq  (\frac{1}{2^{k-1}}-\frac{p}{2})\binom{n-1}{k-1}+O(n^{k-2}).$$
\end{theorem}

We now aim to derive a lower bound for $c(p, n)$ when $p$ is close to $1/2$.
Given two (hyper)graphs $G$ and $H$, denote by $G\cup H$ the disjoint union of $G$ and $H$. For an integer $k\ge2$, denote by $K_n^{(k)}$ the complete $k$-uniform hypergraph on $n$ vertices. In particular, write $K_{n}$ instead of $K_n^{(2)}$ if $k=2$.
\begin{theorem}\label{33}
	Let $G$ be a graph with $n$ vertices and $p\binom{n}{2}$ edges, where $\frac{1}{4}\leq p \leq \frac{1}{2}$. Suppose that
	\[
	G:=
	K_{\sqrt{p} n}\cup \left(\frac{\sqrt p-p}{2}n\right)K_2\cup\left((1-\sqrt p)^2n\right)K_1.
	\]
	Then
	\[
	\operatorname{grdisc}(G)\ge
	\left(\sqrt{p}-p\right)(n-1)+\sqrt{p}-1.
	\]
\end{theorem}
\noindent\textbf{Remark}. For clarity, we ignore all rounding up and down. Indeed, there do exist $p$ and $n$ satisfying this theorem without rounding up and down. In particular, if $p=1/4$ and $n$ is multiple of 8, then our theorem implies that $G:=K_{n/2}\cup (n/8)K_2\cup (n/4)K_1$ and $\operatorname{grdisc}(G)\ge (n-3)/4$.
% (for example, $p=1/4$ and $n=0\;(\text{mod}\; 8)$). %For general, we provide in the next theorem.
\begin{proof}[{\bf Proof of Theorem \ref{33}}]
	We first construct an induced subgraph $G_{n-1}$ by deleting  a vertex $v_n$  from $G$ such that $\operatorname{grdisc}(G_{n-1},G)$ is minimized. Clearly,
	\begin{align*}
		\operatorname{grdisc}(G_{n-1},G)&=\left| p\binom{n}{2}-d(v_n,G)-p\binom{n-1}{2}\right|\\
		&=\left|p(n-1)-d(v_n,G)\right|.
	\end{align*}
	This further implies that
	\begin{align*}
		\operatorname{grdisc}(G)&=\min_{\sigma\in S_n}\max_{i\in [n]} \operatorname{grdisc}(G_{\sigma(i)},G)\\
		&\geq \min_{\sigma(n)=v_n}\operatorname{grdisc}(G_{n-1},G)\\
		%&\geq\min_{v_n} \left| e(G_{n-1})-\frac{1}{4}\binom{n-1}{2}\right|\\
		%&\geq\min_{v_n}\left|\frac{1}{4}\binom{n}{2}-d(v_n,G_n)-\frac{1}{4}\binom{n-1}{2}\right|\\
		&=\min_{v_n\in V(G)}\left|p(n-1)-d(v_n,G)\right|.
		%&=\min\left\{\left|\frac{n-3}{4}\right|, \left|\frac{n-9}{4}\right|, \left|\frac{n-1}{4}\right|\right\}\\ &=  \frac{1}{4}(n-9).
	\end{align*}
	Note that  $d(v_n,G)\in\{\sqrt{p}n-1,1,0\}$. Thus, the desired result holds.
\end{proof}

Now, we give the following general constructions for any $k\ge2$.
\begin{theorem}\label{34}
	For a real $\frac{1}{2^k} \leq p \leq \frac{1}{2}$ and an integer $k\ge2$, let  $s=p^{1/k}\le(1-1/\sqrt[k]{n})^k$. For sufficiently large $n$, suppose that $H$ is a $k$-uniform hypergraph with $n$ vertices and $p\binom{n}{k}$ edges satisfying
	%\[G:=
	%    \begin{cases}
		%		\frac{1}{\sqrt[k-1]{2^kp}}K^{(k)}_{\sqrt[k-1]{2p}n}\bigcup\frac{1-2p}{4}nK_2\bigcup pnK_1, \;\text{if} \,\;p\le1/4,\\
		%		K^{(k)}_{\sqrt[k]{p} n}\bigcup \frac{\sqrt p-p}{2}nK_2\bigcup(1-\sqrt p)^2nK_1, \;\text{if}\,\; p>1/4.
		%	\end{cases}\]
	\[
	H:=
	K^{(k)}_{\lfloor s n\rfloor}\cup H',
	\]
	where $H'$ is a $k$-uniform hypergraph with its maximum degree $\Theta(s^{k-1} n^{k-2})$.
	Then
	\[
	\operatorname{grdisc}(H)\ge
	\left(p^{(k-1)/k}-p-o(p)\right)\binom{n-1}{k-1}.
	\]
\end{theorem}
\noindent\textbf{Remark}. In particular, if $p=2^{-k}$ and $n$ is even, then our theorem implies that $H:=K_{n/2}\cup H'$ for some $k$-uniform hypergraph $H'$ on $n/2$ vertices with its maximum degree $o(n^{k-1})$, and $\operatorname{grdisc}(H)\ge(1/2^k-o(1))\binom{n-1}{k-1}$.
\begin{proof}[{\bf Proof of Theorem \ref{34}}]
	Note first that
	\begin{align*}
		\binom{\lfloor s n\rfloor-1}{k-1}= p^{(k-1)/k}\binom{n-1}{k-1}-\Theta\left(p^{(k-2)/k}n^{k-2}\right).
	\end{align*}
	This implies that
	\begin{align*}
		e\left(K^{(k)}_{\lfloor s n\rfloor}\right)=\frac{\lfloor sn\rfloor}{k}\binom{\lfloor s n\rfloor-1}{k-1}= p\binom{n}{k}-\Theta\left(p^{(k-1)/k}n^{k-1}\right).
	\end{align*}
	Observe that  $|V(H')|=(1-s)n+\Theta(1)$. This implies that $H$ is well-defined by choosing $H'$ properly under the constraint of its maximum degree $\Theta(s^{k-1} n^{k-2})$.
	
	Now, we construct an induced subgraph $H_{n-1}$ by deleting  a vertex $v_n$  from $H$ such that $\operatorname{grdisc}(H_{n-1},H)$ is minimized. Clearly,
	\begin{align*}
		\operatorname{grdisc}(H_{n-1},H)&=\left| p\binom{n}{k}-d(v_n,H)-p\binom{n-1}{k}\right|\\
		&=\left|p\binom{n-1}{k-1}-d(v_n,H)\right|.
	\end{align*}
	This further implies that
	\begin{align*}
		\operatorname{grdisc}(H)&=\min_{\sigma\in S_n}\max_{i\in [n]} \operatorname{grdisc}(H_{\sigma(i)},H)\\
		&\geq \min_{\sigma(n)=v_n}\operatorname{grdisc}(H_{n-1},H)\\
		%&\geq\min_{v_n} \left| e(G_{n-1})-\frac{1}{4}\binom{n-1}{2}\right|\\
		%&\geq\min_{v_n}\left|\frac{1}{4}\binom{n}{2}-d(v_n,G_n)-\frac{1}{4}\binom{n-1}{2}\right|\\
		&=\min_{v_n\in V(G)}\left|p\binom{n-1}{k-1}-d(v_n,H)\right|.
		%&=\min\left\{\left|\frac{n-3}{4}\right|, \left|\frac{n-9}{4}\right|, \left|\frac{n-1}{4}\right|\right\}\\ &=  \frac{1}{4}(n-9).
	\end{align*}
	Note that  $d(v_n,H)=\binom{\lfloor s n\rfloor-1}{k-1}=p^{(k-1)/k}\binom{n-1}{k-1}-\Theta(p^{(k-2)/k}n^{k-2})$ or $d(v_n,H)\le \Delta(H')=\Theta(p^{(k-1)/k}n^{k-2})$. Thus, we have the desired result.
\end{proof}

\section*{Acknowledgements} The authors are grateful to the anonymous reviewers for helpful comments on an earlier version of this manuscript.

%{\bf Acknowledgement.}


\begin{thebibliography}{99}
	\bibitem{BS1995} J. Beck, V. T. S\'os, Discrepancy theory, In Handbook of Combinatorics, Vol. 2, Elsevier, Amsterdam (1995) 1405--1446.
	
	%\bibitem{BJS2017} B. Bollob\'as, S. Janson, A.D. Scott, Packing random graphs and hypergraphs. Random Struct. Alg. 1 (2017) 3--13.
	
	\bibitem{BS2002} B. Bollob\'as, A. Scott, Problems and results on judicious partitions, Random Struct. Alg. 21 (2002) 414--430.
	
	\bibitem{BS2006} B. Bollob\'as, A.D. Scott, Discrepancy in graphs and hypergraphs, in: Ervin Gyori, Gyula O.H. Katona, Laszlo Lov\'asz (Eds.), More Sets, Graphs and Numbers, in: Bolyai Soc. Math. Stud., vol. 15, Springer, Berlin (2006) 33--56.
	
	\bibitem{BS2011} B. Bollob\'{a}s, A. D. Scott, Intersections of graphs, J. Graph Theory 66 (2011) 261--282.
	
	\bibitem{BS2015JCTB} B. Bollob\'as, A.D. Scott, Intersections of hypergraphs, J. Combin. Theory, Series B 110 (2015) 180--208.
	
	\bibitem{BS2015EJC}	B. Bollob\'as, A.D. Scott, Intersections of random hypergraphs and tournaments, European J. Combin. 44 (2015) 125--139.
	
	\bibitem{Cha2000} B. Chazelle, The discrepancy method, Cambridge University Press, Cambridge (2000) pp. xviii+463.
	
	\bibitem{ES1971} P. Erd\H{o}s, J. Spencer, Imbalances in $k$-colorations, Networks 1 (1971/1972) 379--385.
	
	\bibitem{EGP1988} P. Erd\H{o}s, M. Goldberg, J. Pach, J. Spencer, Cutting a graph into two dissimilar halves, J. Graph Theory 12 (1988) 121--131.
	
	%\bibitem{LLB2013} C. Lee, P. Loh, and B. Sudakov, Self-similarity of graphs, SIAM J. Discrete Math. 27 (2013) 959--972.
	
	\bibitem{MNS2013} J. Ma, H. Naves, B. Sudakov, Discrepancy of random graphs and hypergraphs, Random Struct. Alg. 47 (2015) 147--162.
	
	\bibitem{Mat1999} J. Matou\v{s}ek, Geometric discrepancy, Algorithms and Combinatorics, Vol. 18, Springer-Verlag, Berlin (1999) pp. xii+288.
	
	\bibitem{RST2023} E. R\"aty, B. Sudakov, I. Tomon, Positive discrepancy, MaxCut, and eigenvalues of graphs, Transactions of the AMS (2026), to appear.
	
	\bibitem{RT2024} E. R\"aty, I. Tomon, Bisection Width, discrepancy, and eigenvalues of hypergraphs, J. Combin. Theory, Series B 177 (2026) 186-215.
	
	\bibitem{Sos1983} V. T. S\'os, Irregularities of partitions: Ramsey theory, uniform distribution, In Surveys in Combinatorics (Southampton, 1983), London Math. Soc. Lecture Note Ser., 82, Cambridge Univ. Press, Cambridge-New York (1983) 201--246.
\end{thebibliography}
\end{document}